\documentclass[12pt]{amsart}

\usepackage{amsmath,amsthm,amssymb}
\font\teneufm=eufm10
\font\seveneufm=eufm7
\font\fiveeufm=eufm5
\newfam\frakturfam

\textfont\frakturfam=\teneufm
\scriptfont\frakturfam=\seveneufm
\scriptscriptfont\frakturfam=\fiveeufm

\newtheorem{lm}{Lemma}
\newtheorem{theor}{Theorem}
\newtheorem{co}{Corollary}

\def\bee{\begin{eqnarray}}
\def\bes{\begin{eqnarray*}}
\def\eee{\end{eqnarray}}
\def\ees{\end{eqnarray*}}

\def\a{\alpha}
\def\b{\beta}
\def\te{\theta}
\def\g{\gamma}

\def\s{\sigma}
\def\t{\tau}

\def\Proof{{\sl Proof.}\ }

\title{Defining relations for automorphism groups of free algebras}

\begin{document}
\date{}
\maketitle

\begin{center}

{\bf U.\,U.\,Umirbaev}\\
Eurasian National University,\\
 Astana, 010008, Kazakhstan \\
e-mail: {\em umirbaev@yahoo.com}

\end{center}

\begin{abstract}
We describe a set of defining relations for automorphism groups of finitely generated free algebras 
of Nielsen-Schreier varieties. In particular, this gives a representation of the automorphism groups of free Lie algebras 
by generators and defining relations. 
\end{abstract}



\section{Introduction}

\hspace*{\parindent}

It is well known \cite{Czer,Jung,Kulk,Makar} that all automorphisms 
of polynomial algebras and free associative algebras in two variables are tame. Moreover, 
 the groups of automorphisms 
of polynomial algebras and free associative algebras in two variables are isomorphic  
and have a nice representation as a free product of groups (see, for example \cite{Cohn,Essen}). 

It was recently proved    
 that the automorphism groups 
of polynomial algebras \cite{SU4,Umi25,SU3} and free associative algebras \cite{Umi31,Umi33} in three variables over a 
field of characteristic $0$ cannot be generated by all 
 elementary automorphisms, i.e. there exist wild automorphisms.  
 Defining relations of the tame automorphism group 
of polynomial algebra in three variables were described in \cite{Umi31,Umi32}. 

There are several well-known descriptions 
of the automorphism group of a free group by generators and defining relations (see, for example \cite{Magnus}). 
 P.Cohn  proved \cite{Cohn2} that all automorphisms of finitely generated free 
Lie algebras are tame. 
 Later this result was extended to free algebras of Nielsen-Schreier varieties \cite{Lewin}. Recall that 
a variety of universal algebras is called Nielsen-Schreier, if any subalgebra of a free algebra of this variety is free, 
i.e. an analog of the classical Nielsen-Schreier theorem is true. 
 The varieties of all nonassociative algebras \cite{Kurosh}, 
 commutative and anticommutative algebras \cite{Shirshov}, Lie algebras 
 \cite{Shirshov,Witt} are Nielsen-Schreier. Other examples of Nielsen-Schreier 
 varieties can be found in \cite{Mikhalev,SU1,Stern,Umi11}. 

So, the automorphism groups of free algebras of Nielsen-Schreier varieties are generated by all 
elementary automorphisms. In this paper we describe a set of defining relations of these
groups. In fact, we show that the relations for elementary automorphisms 
studied in \cite{Umi31,Umi32} are defining relations in this case.

The paper is organized as follows. 
In Section 2 we describe a set of relations for elementary automorphisms and repeat the proofs of two lemmas from \cite{Umi32} for completeness. 
In Section 3 we give some well-known definitions and theorems about free algebras. 
 In Section 4 we prove the main result.

\section{Defining relations}

\hspace*{\parindent}

Let $F$ be an arbitrary field, and let ${\mathfrak M}$ be an arbitrary variety of linear algebras over $F$. By $A=F_{\mathfrak M}<x_1,x_2,\ldots,x_n>$
 denote the free algebra of ${\mathfrak M}$ with a free set of generators $X=\{x_1,x_2,\ldots,x_n\}$, and  
 by $Aut\,A$ denote the group of all automorphisms of this algebra. Let  $\phi = (f_1,f_2,\ldots,f_n)$ denote an automorphism 
  $\phi$ of $A$ such that $\phi(x_i)=f_i,\, 1\leq i\leq n$. An automorphism  
\bee\label{af1} 
\s(i,\a,f) = (x_1,\ldots,x_{i-1}, \a x_i+f, x_{i+1},\ldots,x_n),
\eee
where $0\neq\a\in F,\ f\in F_{\mathfrak M}<X\setminus \{x_i\}>$, 
 is called {\em elementary}.    The subgroup $TA(A)$ of $Aut\,A$ generated by all 
  elementary automorphisms is called the {\em tame automorphism group}, 
 and the elements of this subgroup are called {\em tame automorphisms} 
 of $A$. Nontame automorphisms of $A$ are called {\em wild}.

 Now we describe some relations for elementary automorphisms (\ref{af1}). It is easy to check that   
\bee\label{af2} 
\s(i,\a,f) \s(i,\b,g) = \s(i,\a \b, \b f+g). 
\eee

Note that from this we obtain trivial relations  $\s(i,1,0) = id$, where $1\leq i\leq n$.

If $i\neq j$ and $f\in F_{\mathfrak M}<X\setminus \{x_i,x_j\}>$, then we have also 
\bee\label{af3} 
\s(i,\a,f)^{-1} \s(j,\b,g) \s(i,\a,f) = \s(j,\b,\s(i,\a,f)^{-1}(g)). 
\eee

Consequently, if $i\neq j$ and $f,g\in F_{\mathfrak M}<X\setminus \{x_i,x_j\}>$, then  the automorphisms $\s(i,\a,f)$,  
$\s(j,\b,g)$ commute.

For every pair of integers $k,s$, where $1\leq k\neq s\leq n$, we define a tame automorphism $(ks)$ by putting  
\bes
(ks)=\s(s,-1,x_k) \s(k,1,-x_s) \s(s,1,x_k). 
\ees 
Note that the automorphism $(ks)$ of the algebra $A$ just permutes the variables $x_k$ and $x_s$.  
Now it is easy to see that  
\bee\label{af4}
\s(i,\a,f)^{(ks)}=\s(j,\a,(ks)(f)),
\eee
where $x_j=(ks)(x_i)$. 

\smallskip 

Let $G(A)$ be the abstract group with generators (\ref{af1}) and defining relations (\ref{af2})--(\ref{af4}). 

\begin{lm}\label{al1}
The subgroup of $G(A)$ generated by all 
elements $(ks)$, where $1\leq k\neq s\leq n$, 
is isomorphic to the symmetric group  $S_n$.  
\end{lm}
\Proof
By (\ref{af2}) and (\ref{af3}), we have 
\bes
(ks)^2=\s(s,-1,x_k) \s(k,1,-x_s) \s(s,1,x_k) \s(s,-1,x_k) \s(k,1,-x_s) \s(s,1,x_k) \\
= \s(s,-1,x_k) \s(k,1,-x_s) \s(s,-1,0) \s(k,1,-x_s) \s(s,1,x_k) \\
= \s(s,-1,x_k) \s(s,-1,0) \s(k,1,-x_s)^{\s(s,-1,0)} \s(k,1,-x_s) \s(s,1,x_k) \\
= \s(s,1,-x_k) \s(k,1,x_s) \s(k,1,-x_s) \s(s,1,x_k) 
= \s(s,1,-x_k) \s(s,1,x_k) =id. 
\ees
Then (\ref{af4}) gives 
\bes
(ks)^{(sk)}=\s(s,-1,x_k)^{(sk)} \s(k,1,-x_s)^{(sk)} \s(s,1,x_k)^{(sk)} \\
 =
\s(k,-1,x_s) \s(s,1,-x_k) \s(k,1,x_s)=(sk), 
\ees 
i.e. $(ks)=(sk)$. Now it is not difficult to deduce from (\ref{af2})--(\ref{af4}) that  
\bes
[(ij),(ks)]=id, \, \, (ik)^{(is)}=(ks), 
\ees
where $i,j,k,s$ are all distinct. It is 
immediate that the given relations imply the defining relations of the group $S_n$ 
with respect to the system of generators $(i\,\,i+1)$, where $1\leq i\leq n-1$, which are indicated in \cite{CM}. 
$\Box$

 By Lemma  \ref{al1}, the elements of the symmetric group $S_n$ can be identified with  elements of $G(A)$.  
Note that (\ref{af4}) can be rewritten as   
\bes
\s(i,\a,f)^{\pi}=\s(\pi^{-1}(i),\a,\pi^{-1}(f)),
\ees
where $\pi \in S_n$.

\smallskip 

It is well known that the group of affine automorphisms $Af_n(F)$ of the algebra $A$ is generated by all affine elementary automorphisms.  
\begin{lm}\label{al2}
The relations (\ref{af2})--(\ref{af4}) for elementary affine automorphisms are defining relations of the group $Af_n(F)$. 
\end{lm}
\Proof
Let $\varphi$ be a product of elementary affine automorphisms.  
Suppose that $\varphi=id$. 
 By (\ref{af2}) and (\ref{af3}), we can represent $\varphi$ in the form  
\bes
\varphi=\s(1,1,\a_1)\s(2,1,\a_2)\ldots \s(n,1,\a_n) \varphi', 
\ees
where $\varphi'$ is a product of elementary linear automorphisms. Obviously, $\a_1=\a_2=\ldots =\a_n=0$. 
Therefore we can assume that $\varphi$ is a product of elementary linear automorphisms.  
 By (\ref{af2}) and (\ref{af3}), we can easily represent  $\varphi$ in the form  
\bes
\varphi=\s(1,\a_1,0)\s(2,\a_2,0)\ldots \s(n,\a_n,0) \varphi', 
\ees
where $\varphi'$ is a product of elementary automorphisms of the type $\s(i,1,f)$. By (\ref{af2})--(\ref{af4}), we have   
\bes
\s(k,\a,0)=\s(s,\a,0)^{(ks)} \\
=\s(s,-1,x_k) \s(k,1,-x_s) \s(s,1,x_k) \s(s,\a,0) \s(s,-1,x_k) \s(k,1,-x_s) \s(s,1,x_k) \\
=\s(s,-1,x_k) \s(k,1,-x_s) \s(s,-\a,0) \s(s,1,(1-\a)x_k) \s(k,1,-x_s) \s(s,1,x_k) \\
=\s(s,-1,x_k) \s(s,-\a,0) \s(k,1,\a^{-1}x_s) \s(s,1,(1-\a)x_k) \s(k,1,-x_s) \s(s,1,x_k) \\
=\s(s,\a,0) \s(s,1,-\a x_k) \s(k,1,\a^{-1}x_s) \s(s,1,(1-\a)x_k) \s(k,1,-x_s) \s(s,1,x_k).
\ees
By using this relation, we can represent $\varphi$ in the form   
\bes
\varphi=\s(n,\b_n,0) \varphi', 
\ees
where $\varphi'$ is a product of elementary linear automorphisms of the form $\s(i,1,f)$. Hence $\b_n=1$. 
Note that $\s(i,1,f)$ can be represented as a product of automorphisms 
\bee\label{af5}
X_{ij}(\lambda)=\s(j,1,\lambda x_i), \ \ \lambda \in F, \ \ i\neq j. 
\eee
Thus, we can assume that $\varphi$ is a product of automorphisms of the form (\ref{af5}). 

Let $G$ be the subgroup of $TA(A)$ generated by all automorphisms of the form (\ref{af5}). 
We define a map 
\bes
J : G \longrightarrow SL_n(F), 
\ees
where $J(\psi)$ is the Jacobian matrix of $\psi\in G$. By $e_{ij}$ denote the standard matrix units and by  
 $E_{ij}(\lambda)$ denote the elementary matrix $E+\lambda e_{ij}$, where $E$ is the unit matrix, $i \neq j$, and $\lambda \in F$. 
It is easy to check that 
\bes
J(X_{ij}(\lambda))=E_{ij}(\lambda), 
\ees 
and that $J$ is an isomorphism of groups. 

Now it is sufficient to prove that every relation of the group $SL_n(F)$ is a corollary of (\ref{af2})--(\ref{af4}). 
Obviously,  (\ref{af2})--(\ref{af3}) cover the Steinberg relations (see, for example \cite{Milnor}). 
Besides, according to \cite{Milnor}, we need to check the relation  
\bes
\{u,v\}=id, \ \ 0\neq u,v\in F, 
\ees 
where 
\bes
\{u,v\}=h_{ij}(uv) h_{ij}(u)^{-1} h_{ij}(v)^{-1}, \\
h_{ij}(u)=w_{ij}(u) w_{ij}(-1), \\
w_{ij}(u)=X_{ij}(u) X_{ji}(-u^{-1}) X_{ij}(u). 
\ees
Applying (\ref{af2})--(\ref{af4}) we have 
\bes
w_{ij}(u)= \s(j,1,u x_i) \s(i,1,-u^{-1} x_j) \s(j,1,u x_i) \\
= \s(j,1,u x_i) \s(i,u,0) \s(i,1,-x_j) \s(i,u^{-1},0) \s(j,1,u x_i) \\
= \s(i,u,0) \s(j,1,u x_i)^{\s(i,u,0)} \s(i,1,-x_j) \s(j,1,u x_i)^{\s(i,u,0)} \s(i,u^{-1},0) \\
= \s(i,u,0) \s(j,1,x_i) \s(i,1,-x_j) \s(j,1,x_i) \s(i,u^{-1},0) \\
= \s(i,u,0) \s(j,-1,0) \s(j,-1,x_i) \s(i,1,-x_j) \s(j,1,x_i) \s(i,u^{-1},0) \\
= \s(i,u,0) \s(j,-1,0) (ij) \s(i,u^{-1},0) = (ij) \s(i,u,0)^{(ij)} \s(j,-1,0)^{(ij)} \s(i,u^{-1},0)\\
= (ij) \s(j,u,0) \s(i,-1,0) \s(i,u^{-1},0)= (ij) \s(j,u,0) \s(i,-u^{-1},0). 
\ees
Consequently, 
\bes
h_{ij}(u)=w_{ij}(u) w_{ij}(-1) = (ij) \s(j,u,0) \s(i,-u^{-1},0) (ij) \s(j,-1,0) \s(i,1,0) \\
= \s(j,u,0)^{(ij)} \s(i,-u^{-1},0)^{(ij)} \s(j,-1,0)  \\ 
= \s(i,u,0) \s(j,-u^{-1},0) \s(j,-1,0) = \s(i,u,0) \s(j,u^{-1},0).  
\ees
Hence  
\bes
\{u,v\}=h_{ij}(uv) h_{ij}(u)^{-1} h_{ij}(v)^{-1} \\
=  \s(i,uv,0) \s(j,(uv)^{-1},0) \s(i,u,0) \s(j,u^{-1},0) \s(i,v,0) \s(j,v^{-1},0)=id.
\ees
 Thus we can say that every relation of the group $SL_n(F)$ follows from (\ref{af2})--(\ref{af4}). 
$\Box$

\section{Reductions of automorphisms}

\hspace*{\parindent}

Let ${\mathfrak M}$ be an arbitrary homogeneous variety of linear algebras over a field $F$. 
Recall that if $F$ is infinite, then any variety of linear algebras over
$F$ is homogeneous \cite{KBKA}.  Let $A=F_{\mathfrak M}<x_1,x_2,\ldots,x_n>$
 be the free algebra of ${\mathfrak M}$ with free generators $x_1,x_2,\ldots,x_n$. 
The highest homogeneous part  $\overline{f}$ and the degree  $\deg f$ can be defined in the usual way. 
If $f_1,f_2,\ldots,f_k \in A$, then denote by 
$<f_1,f_2,\ldots,f_k>$ the subalgebra of $A$ generated by these elements. 
 
Let $\te =(f_1,f_2,\ldots,f_k)$ be an arbitrary $k$-tuple of elements of the algebra  $A$. The number  
\bes
\deg \te=\deg f_1+\deg f_2+\ldots+\deg f_k
\ees
is called the {\em degree} of $\te$.

Recall that an {\em elementary transformation} of a $k$-tuple   $\te =(f_1,f_2,\ldots,f_k)$
 is, by definition, a transformation that changes only one element $f_i$ to an element of the form $\a
f_i+g$, where  $0\neq \a\in F$ and $g\in \langle\{f_j|j\neq
i\}\rangle$. The notation  
\bes
\te \rightarrow \tau
\ees 
means that the $k$-tuple $\tau$ is obtained from $\te$ by a single elementary transformation. 
A $k$-tuple $\te$ is called {\em  elementarily reducible} or {\em admits an elementary reduction} if there exists a $k$-tuple
     $\tau$ such that $\te\rightarrow \tau$ and $\deg\tau<\deg\te$. The element  $f_i$ of the $k$-tuple $\te$ which was changed in $\tau$   
 to an element of less degree is called  {\em reducible} 
 and we will say also that  $f_i$ is {\em reduced in $\te$ by the $k$-tuple $\tau$}.

Consider a finite number of elements 
\bee\label{af6} 
f_1,f_2,\ldots,f_k
\eee
of the algebra $A$. The elements (\ref{af6}) are called {\em free} if the 
subalgebra $<f_1,f_2,\ldots,f_k>$ of $A$ is a free algebra of the 
variety ${\mathfrak M}$ with free generators (\ref{af6}). If 
\bes
\overline{f_i} \notin <\overline{f_1},\ldots,\overline{f_{i-1}},\overline{f_{i+1}},\ldots,\overline{f_n}>
\ees
for any $i$, then the elements (\ref{af6}) are called {\em reduced}. 

From any $k$-tuple $(f_1,f_2,\ldots,f_k)$ by several elementary transformations 
we can get a $k$-tuple $(g_1,g_2,\ldots,g_s,0,\ldots,0)$, where $s\leq k$, 
such that the elements $g_1,g_2,\ldots,g_s$ are reduced. Note that $<f_1,f_2,\ldots,f_k>=<g_1,g_2,\ldots,g_s>$. 

The statement of the next lemma is well known (see, 
for example \cite{Shirshov}) and very useful in studying free algebras.

\begin{lm}\label{al3}
Assume that the elements $\overline{f_1},\overline{f_2},\ldots,\overline{f_k}$ are free. 
If $f\in <f_1,f_2,\ldots,f_k>$, then $\overline{f}\in <\overline{f_1},\overline{f_2},\ldots,\overline{f_k}>$. 
\end{lm}

From now on we assume that ${\mathfrak M}$ is a homogeneous Nielsen-Schreier variety of linear algebras. 
The main property of Nielsen-Schreier varieties is given in the next lemma (see, for example \cite{Lewin}). 
\begin{lm}\label{al4}
Assume that $f_1,f_2,\ldots,f_k$ are homogeneous elements of $A$ and $\deg\,f_1\leq \deg\,f_2,\leq \ldots \leq \deg\,f_k$. 
If the elements $f_1,f_2,\ldots,f_k$ are not free, then there exists $i$ such that 
$f_i\in <f_1,f_2,\ldots,f_{i-1}>$. 
\end{lm}

\begin{co}\label{cn1}
Any finite reduced system of elements of the algebra $A$ is free. 
\end{co}

Note that the statement of this corollary for infinite systems of elements is also true \cite{Lewin}. 
Free systems of elements in free algebras were studied in \cite{MShZ,Umi4} via Fox derivatives. 

\begin{co}\label{cn2}
Any automorphism of the algebra $A$ of degree more than $n$ is elementarily reducible. 
\end{co}

\begin{co}\label{cn3}
Automorphisms of the algebra $A$ are tame. 
\end{co}

Suppose that $\te =(f_1,f_2,\ldots,f_n)$ and $\s(i,\a,f)$ is an elementary automorphism of the form (\ref{af1}). If 
\bes
\tau = (f_1,\ldots,f_{i-1},\a f_i+f(f_1,\ldots,f_n), f_{i+1},\ldots,f_n),
\ees
then instead of $\te\rightarrow \tau$ we often write  
\bes
\te\stackrel{\s(i,\a,f)}\longrightarrow \tau. 
\ees

Assume that   
\bee\label{af7}
\te=\phi_1\phi_2\ldots \phi_r \in Aut\,A, 
\eee
where $\phi_i$, $1\leq i\leq r$, are elementary automorphisms. 
Put  
\bes
\psi_i=\phi_1\phi_2\ldots\phi_i, \,\, 0\leq i\leq r.  
\ees
In particular, we have 
\bes
\psi_r=\te,\,\,\psi_0=id.  
\ees
To (\ref{af7}) corresponds the sequence of elementary transformations  
\bee\label{af8}
id=\psi_0\stackrel{\phi_1}{\rightarrow}\psi_1\stackrel{\phi_2}{\rightarrow}\psi_2\stackrel{\phi_3}{\rightarrow}\ldots \stackrel{\phi_r}{\rightarrow}\psi_r=\te.  
\eee
So, every tame automorphism  $\te$ has a sequence of elementary transformations  of the form (\ref{af8}). 
If $\deg\,\te > n$ and $\deg\,\psi_i < \deg\,\te$ for any $i<r$, 
then the sequence (\ref{af8}) will be called a {\em minimal representation} of $\te$. 
Note that the representations (\ref{af7}) and (\ref{af8}) of the automorphism  $\te$ are equivalent. 
If (\ref{af8}) is a minimal representation of $\te$, then the representation (\ref{af7}) will be also called a {\em minimal representation} of $\te$.

\section{The main result}

\hspace*{\parindent}

As above, ${\mathfrak M}$ is a homogeneous Nielsen-Schreier variety of linear algebras over a field $F$, and $A=F_{\mathfrak M}<x_1,x_2,\ldots,x_n>$
 is a free algebra of ${\mathfrak M}$. We know that $TA(A)=Aut\,A$ and 
the elementary automorphisms (\ref{af1}) are generators of the group $Aut\,A$. 

\begin{theor}\label{at1} 
The relations (\ref{af2})--(\ref{af4}) are defining relations of the group $Aut\,A$ 
with respect to the generators (\ref{af1}). 
\end{theor}
{\bf Beginning of the proof.}
Assume that   
\bee\label{af9}
\varphi_1\varphi_2\ldots\varphi_k=id=(x_1,x_2,\ldots,x_n),
\eee
where $\varphi_i$, $1\leq i\leq k$, are elementary automorphisms. Put    
\bes
\te_i=\varphi_1\varphi_2\ldots\varphi_i,\,\,0\leq i\leq k.   
\ees
In particular, we have $\te_0=\te_k=(x_1,x_2,\ldots,x_n)$. To (\ref{af9}) corresponds the sequence of elementary transformations  
\bee\label{af10}
id=\te_0\stackrel{\varphi_1}{\rightarrow} \te_1\stackrel{\varphi_2}{\rightarrow} \ldots \stackrel{\varphi_k}{\rightarrow} \te_k=id.
\eee

Put $d=max\{\deg\,\te_i | 0\leq i\leq k\}$. Let $i_1$ be the minimal number and $i_2$ be the maximal 
number which satisfy the equations $\deg\,\te_{i_1}=d$ and $\deg\,\te_{i_2}=d$. Put $q=i_2-i_1$. 
The pair $(d,q)$ will be called the {\em exponent} of the relation (\ref{af9}). 

To prove the theorem, we show that (\ref{af9}) 
 follows from (\ref{af2})--(\ref{af4}). Assume that our theorem is not true. 
 Call a relation of the form (\ref{af9})  
  {\em trivial} if it follows from (\ref{af2})--(\ref{af4}). 
 We choose a nontrivial relation (\ref{af9})  with the minimal exponent $(d,q)$ 
 with respect to the lexicographic order. To arrive at a contradiction, 
 we show that (\ref{af9}) is also trivial.  

If $d=n$, then Lemma \ref{al2} gives the triviality of the relation (\ref{af9}). 
Therefore we can assume that $d>n$.

Our plan is to change the product  (\ref{af9}) by using (\ref{af2})--(\ref{af4}) and 
to obtain a new sequence (\ref{af10}) whose exponent is strictly less than $(d,q)$. 
 Below we prove Lemmas \ref{al5}--\ref{al14} and then complete the proof of the theorem. $\Box$

\medskip

 Denote by $t=[\frac{q}{2}]$ the integral part of $\frac{q}{2}$.
Put also   
\bes
\phi=\te_{i_1+t-1},\,\, 
\te=\te_{i_1+t},\,\, \tau=\te_{i_1+t+1}, \,\, \s_1=\varphi_{i_1+t}, \,\, \s_2=\varphi_{i_1+t+1}. 
\ees
 Then we have  
\bee\label{af11}
\phi \stackrel{\s_1}{\longrightarrow} \theta \stackrel{\s_2}{\longrightarrow} \tau. 
\eee

\begin{lm}\label{al5} 
The following statements are true:  
\begin{itemize}
\item[(1)]  $d=\deg\,\te$, $t=0$, and  
\bee\label{af12}
\te=\varphi_1\varphi_2\ldots \varphi_{i_1+t} 
\eee  
is a minimal representation of $\te$. 
\item[(2)]  If $q=0$, then  
\bes
\te=\varphi_{k}^{-1}\varphi_{k-1}^{-1}\ldots \varphi_{i_1+t+1}^{-1}  
\ees
is also a minimal representation of $\te$. 
\item[(3)]  If $q=1$, then $(d(\tau),t(\tau))=(d,t)$ and  
\bes
\tau=\varphi_{k}^{-1}\varphi_{k-1}^{-1}\ldots \varphi_{i_1+t+2}^{-1}  
\ees
is a minimal representation of $\tau$. Moreover, in (\ref{af9}) 
the product (\ref{af12}) can be replaced by an arbitrary minimal representation of 
 $\te$. 
\end{itemize} 
\end{lm}
\Proof
Assume that $(d(\te),t(\te))<(d,t)$ and let (\ref{af7}) be a minimal representation of $\te$. Then (\ref{af9}) is a consequence of the equalities    
\bee\label{af13}
\varphi_1\varphi_2\ldots \varphi_{i_1+t}\phi_r^{-1}\ldots \phi_2^{-1}\phi_1^{-1} =id,  
\eee  
\bee\label{af14}
\phi_1\phi_2\ldots\phi_r\varphi_{i_1+t+1}\ldots \varphi_{k-1}\varphi_k=id.  
\eee  
To (\ref{af13}) corresponds the sequence of elementary transformations  
\bes
(x_1,x_2,\ldots,x_n)\rightarrow \te_1\rightarrow\ldots\rightarrow\te_{i_1+t}=
\te=\psi_r\rightarrow\psi_{r-1}\ldots\rightarrow\psi_1\rightarrow(x_1,x_2,\ldots,x_n), 
\ees
and to (\ref{af14}) corresponds  
\bes
(x_1,x_2,\ldots,x_n)\rightarrow\psi_1\rightarrow\ldots\rightarrow\psi_r=\te=\te_{i_1+t}\rightarrow\te_{i_1+t+1}
\rightarrow\ldots\rightarrow\te_{k-1}\rightarrow(x_1,x_2,\ldots,x_n). 
\ees
Since $(d(\te),t(\te))<(d,t)$, it follows that (\ref{af13}) and (\ref{af14}) 
have exponents strictly less than $(d,q)$. This gives the first statement of the lemma. 

It is obvious that  the relation 
\bes
\varphi_k^{-1}\varphi_{k-1}^{-1}\ldots\varphi_1^{-1}=id 
\ees
has the same exponent $(d,q)$. Applying the first statement of the lemma  to this relation, 
we get the second statement of the lemma, as well as the minimality of the representation of $\tau$  
if $q=1$. If $q=1$, then (\ref{af13}) 
has exponent strictly less than $(d,q)$, and  (\ref{af14}) has the exponent $(d,q)$. 
Consequently, (\ref{af9}) and (\ref{af14}) are equivalent modulo  
(\ref{af2})--(\ref{af4}). Thus $\te$ can be changed by 
an arbitrary minimal representation in (\ref{af14}). $\Box$

Put $\te=(f_1,f_2,\ldots,f_n)$. According to Lemma \ref{al5},    
 $t=0$, $q=0, 1$, and 
\bes
\deg\,\phi < \deg\,\te = d \geq \deg\,\tau. 
\ees
Without loss of generality, we can assume that 

\bee\label{af15}
\t=(f_1,f_2,\ldots,f_{n-1},f), 
\eee
where 
\bes
\ \  f=\b f_n+B,\ \  B=b(f_1,f_2\ldots,f_{n-1}), \ \ \deg\,B \leq \deg f_n.
\ees

\begin{lm}\label{al6}
If $\phi$ reduces the element $f_n$ of $\te$, then the relation (\ref{af9}) is trivial.  
\end{lm}
\Proof
Applying (\ref{af2}) we can replace $\s_1 \s_2$ 
 by an elementary automorphism.  Obviously, this replacement also decreases the exponent of (\ref{af10}).
$\Box$

By Lemma \ref{al6}, we can assume that $\phi$ reduces one of the elements  $f_1,f_2,\ldots,f_{n-1}$ of $\te$. 
 
\begin{lm}\label{al7} Assume that  $\phi$ reduces the element  $f_i$ of $\te$, where $1\leq i\leq n-1$.  
If $\phi'$ also reduces the element $f_i$ of $\te$, then in 
 (\ref{af11}) the automorphism $\phi$ can be replaced by $\phi'$. 
\end{lm}
\Proof
According to (\ref{af2}), in this case the elementary 
transformation $\phi\rightarrow\te$ can be changed to 
 $\phi\rightarrow\phi'\rightarrow\te$. Since $\deg\,\phi'<\deg\,\te=d$, 
 the exponent $(d,q)$ of the sequence (\ref{af10}) does not change after this replacement. 
 But in the new sequence  (\ref{af10}) we have $\phi'$  instead of $\phi$. $\Box$
 
 From now on we assume that  $\phi$ reduces the element  $f_i$ of $\te$, where $1\leq i\leq n-1$. 
Taking  Lemma \ref{al7} into account, we can assume that  
\bee\label{af16}
\phi=(f_1,\ldots,f_{i-1},g_i,f_{i+1},\ldots,f_n),
\eee
where 
\bes
\ \ g_{i}=f_{i}-C, \ \ C=c(f_1,\ldots,f_{i-1},f_{i+1},\ldots,f_n),\ \  \deg g_{i}<\deg f_{i}.
\ees
Thus, we defined the members of the sequence (\ref{af11}) and we have 
\bes
\s_1=\s(i,1,c(x_1,\ldots,x_{i-1},x_{i+1},\ldots,x_n)), \ \  \s_2=\s(n,\b,b(x_1\ldots,x_{n-1})). 
\ees

\begin{lm}\label{al8}
If the elements $\overline{f_1},\ldots,\overline{f_{i-1}},\overline{f_{i+1}},\ldots,\overline{f_{n-1}}$ are not free, then the relation (\ref{af9}) is trivial. 
\end{lm}
\Proof
If the elements $\overline{f_1},\ldots,\overline{f_{i-1}},\overline{f_{i+1}},\ldots,\overline{f_{n-1}}$ are not free, then according to Lemma \ref{al4}, 
there exists $r$ such that 
\bes
\overline{f_r} \in <\overline{f_1},\ldots,\overline{f_{r-1}},\overline{f_{r+1}},\ldots,\overline{f_{i-1}},\overline{f_{i+1}},\ldots,\overline{f_{n-1}}>.
\ees 
Suppose that 
\bes
\overline{f_r} = T(\overline{f_1},\ldots,\overline{f_{r-1}},\overline{f_{r+1}},\ldots,\overline{f_{i-1}},\overline{f_{i+1}},\ldots,\overline{f_{n-1}}) 
\ees
 and put 
\bes
 g_r=f_r-T(f_1,\ldots,f_{r-1},f_{r+1},\ldots,f_{i-1},f_{i+1},\ldots,f_{n-1}). 
\ees
Then $\deg\,g_r< \deg\,f_r$. Put also 
\bes
\psi_1=(f_1,\ldots,f_{r-1},g_r,f_{r+1},\ldots,f_{i-1},g_i,f_{i+1},\ldots,f_{n-1},f_n), \\
\psi_2=(f_1,\ldots,f_{r-1},g_r,f_{r+1},\ldots,f_{i-1},f_i,f_{i+1},\ldots,f_{n-1},f_n), \\
\psi_3=(f_1,\ldots,f_{r-1},g_r,f_{r+1},\ldots,f_{i-1},f_i,f_{i+1},\ldots,f_{n-1},f). 
\ees
Then we have the sequence of elementary transformations 
\bee\label{af17}
\phi \stackrel{\delta_1}{\longrightarrow} \psi_1 \stackrel{\delta_2}{\longrightarrow} \psi_2 
\stackrel{\delta_3}{\longrightarrow}\psi_3 \stackrel{\delta_4}{\longrightarrow}\tau,  
\eee
where 
\bes
\delta_4=\s(r,1,T(x_1,\ldots,x_{r-1},x_{r+1},\ldots,x_{i-1},x_{i+1},\ldots,x_{n-1})), \ \ \delta_1=\delta_4^{-1}, \\  
\delta_2=\s(i,1,c(x_1,\ldots,x_{r-1},\delta_4(x_{r}),x_{r+1},\ldots,x_{i-1},x_{i+1},\ldots,x_n)), \\ 
 \delta_3=\s(n,\b,b(x_1,\ldots,x_{r-1},\delta_4(x_{r}),x_{r+1},\ldots,x_{n-1})). 
\ees

By (\ref{af3}) we have 
\bes
\delta_1\delta_2\delta_3\delta_4=\delta_2^{\delta_4}\delta_3^{\delta_4}= \s_1\s_2. 
\ees
Then, we can replace the subsequence (\ref{af11}) of (\ref{af10}) by (\ref{af17}). 
Since $\deg\,\psi_1, \deg\,\psi_2, \deg\,\psi_3 < d$, the new sequence (\ref{af10}) has the 
exponent less than $(d,q)$. Consequently, the relation (\ref{af9}) is trivial. 
$\Box$

\begin{lm}\label{al9}
If $\overline{f_i}\in \langle\overline{f_1},\ldots,\overline{f_{i-1}},\overline{f_{i+1}},\ldots,\overline{f_{n-1}}\rangle$, then (\ref{af9}) is trivial. 
\end{lm}
\Proof
Assume that 
\bes
\overline{f_i}=T(\overline{f_1},\ldots,\overline{f_{i-1}},\overline{f_{i+1}},\ldots,\overline{f_{n-1}}).
\ees
 According to Lemma \ref{al7}, we can suppose that 
 \bes
 g_i=f_i-T(f_1,\ldots,f_{i-1},f_{i+1},\ldots,f_{n-1}).
 \ees
  Then  $\s_1=\s(i,1,T)$.  By (\ref{af3}), we have $\s_1\s_2=\s(n,\b,b_1) \s_1$, 
  where $b_1=\s_1(b)\in F_{\mathfrak M}<x_1,\ldots,x_{n-1}>$. 
After the corresponding replacement in (\ref{af9}), 
 $\te$ is replaced by 
\bes
\te'=(f_1,\ldots,f_{i-1},g_i,f_{i+1},\ldots,f_{n-1},f)
\ees
 in (\ref{af10}). 
Since  $\deg\,\te'<d$, the exponent of (\ref{af9}) is decreased.  $\Box$

\begin{lm}\label{al10}
If $\overline{f_n}\in \langle\overline{f_1},\ldots,\overline{f_{i-1}},\overline{f_{i+1}},\ldots,\overline{f_{n-1}}\rangle$, then the relation (\ref{af9}) is trivial.
\end{lm}
\Proof 
Assume that 
\bes
\overline{f_n}=T(\overline{f_1},\ldots,\overline{f_{i-1}},\overline{f_{i+1}},\ldots,\overline{f_{n-1}})
\ees
 and put 
 \bes
 g_n=f_n-T(f_1,\ldots,f_{i-1},f_{i+1},\ldots,f_{n-1}). 
 \ees
According to (\ref{af3}), we have  
\bes
\s_1=\s(i,1,c(x_1,\ldots,x_{i-1},x_{i+1},\ldots,x_n))= \delta_1 \delta_2 \delta_3, 
\ees
where 
\bes
\delta_1=\s(n,1,-T(x_1,\ldots,x_{i-1},x_{i+1},\ldots,x_{n-1})), \\ 
\delta_2=\s(i,1,c_1(x_1,\ldots,x_{i-1},x_{i+1},\ldots,x_n)),  \\
\delta_3=\s(n,1,T(x_1,\ldots,x_{i-1},x_{i+1},\ldots,x_{n-1})). 
\ees
After the corresponding replacement in (\ref{af9}), the elementary transformation  
$\phi\rightarrow\te$ is replaced by the sequence of elementary transformations  
\bes
\phi \rightarrow \psi_1 \rightarrow \psi_2 \rightarrow \te, 
\ees
where 
\bes
\psi_1=(f_1,\ldots,f_{i-1},g_i,f_{i+1},\ldots,f_{n-1},g_n), \\
\psi_2=(f_1,\ldots,f_{i-1},f_i,f_{i+1},\ldots,f_{n-1},g_n). 
\ees
Since $\deg\,\psi_1, \deg\,\psi_2<d=\deg\,\te$, 
 the new sequence (\ref{af10}) has the same exponent $(d,q)$. 
However, instead of $\phi$ we have $\psi_2$, which reduces the element $f_n$ of $\te$. 
By Lemma \ref{al6}, we obtain the triviality of (\ref{af9}). $\Box$

\begin{lm}\label{al11}
If $B$ does not depend on $f_i$, then (\ref{af9}) is trivial.
\end{lm}
\Proof
It means that $b$ does not depend on $x_i$. By (\ref{af3}) we have  
\bes
\s_1 \s_2=\s(i,1,c)\s(n,\b,b)=\s(n,\b,b)\s(2,1,c_1), 
\ees
where $c_1=\s(n,\b,b)^{-1}(c) \in F_{\mathfrak M}<X\setminus \{x_i\}>$. 
After the corresponding replacement in (\ref{af9}), instead of $\te$ we obtain 
\bes
\psi=(f_1,\ldots,f_{i-1},g_i,f_{i+1},\ldots,f_{n-1},f).
\ees 
Since $\deg\,\psi<d$, this replacement also decreases the exponent of (\ref{af9}). $\Box$

\begin{lm}\label{al12}
If $\overline{f_i}=\g \overline{f_n}+T(\overline{f_1},\ldots,\overline{f_{i-1}},\overline{f_{i+1}},\ldots,\overline{f_{n-1}})$, then (\ref{af9}) is trivial.
\end{lm}
\Proof
By Lemma \ref{al9}, we can assume that $\g\neq 0$. By Lemma \ref{al7}, we can also 
assume that $C=\g f_n+T(f_1,\ldots,f_{i-1},f_{i+1},\ldots,f_{n-1})$. Consequently,  
\bes
g_i&=&f_i-\g f_n-T(f_1,\ldots,f_{i-1},f_{i+1},\ldots,f_{n-1}),\\
f_i&=&g_i+\g f_n+T(f_1,\ldots,f_{i-1},f_{i+1},\ldots,f_{n-1}),\\
f_n&=&-\frac{1}{\g} g_i+\frac{1}{\g}f_i-\frac{1}{\g} T(f_1,\ldots,f_{i-1},f_{i+1},\ldots,f_{n-1}). 
\ees
These equalities give rise to the sequence of elementary transformations 
\bes
\psi_1 \rightarrow \psi_2 \rightarrow \te, 
\ees
where 
\bes
\psi_1= (f_1,\ldots,f_{i-1},f_n,f_{i+1},\ldots,f_{n-1},g_i), \\ 
\psi_2= (f_1,\ldots,f_{i-1},f_i,f_{i+1},\ldots,f_{n-1},g_i).
\ees
We have  
\bes
\s_1=\s(i,1,\g x_n+T)=\s(i,1,\g x_n+T(x_1,\ldots,x_{i-1},x_{i+1},\ldots,x_{n-1})). 
\ees
Applying  (\ref{af2}) and (\ref{af3}) we get  
\bes
\s_1= \s(i,1,\g x_n)\s(i,1,T) 
=\s(i,1,\g x_n)\s(n,-\g,x_i)\s(n,-\frac{1}{\g},\frac{1}{\g}x_i)\s(i,1,T)\\
=\s(i,1,\g x_n)\s(n,-\g,x_i)\s(i,1,T)\s(n,-\frac{1}{\g},\frac{1}{\g}x_i)^{\s(i,1,T)}\\
=\s(i,1,\g x_n)\s(n,-\g,x_i)\s(i,1,T)\s(n,-\frac{1}{\g},\frac{1}{\g}(x_i-T))\\
=\s(i,1,\g x_n)\s(n,-\g,x_i)\s(i,\frac{1}{\g},-\frac{1}{\g}x_n)\s(i,\g,x_n+T)\s(n,-\frac{1}{\g},\frac{1}{\g}(x_i-T)).
\ees
Since the transposition $(in)\in S_n$ can be factored as a product of linear elementary automorphisms   
\bes
(in)=\s(i,1,\g x_n)\s(n,-\g,x_i)\s(i,\frac{1}{\g},-\frac{1}{\g}x_n), 
\ees
we obtain   
\bes
\s_1=(in)\s(i,\g,x_n+T)\s(n,-\frac{1}{\g},\frac{1}{\g}(x_i-T)).
\ees
Then  
\bee\label{af18}
\te=(in)\varphi_1^{(in)}\varphi_2^{(in)}\ldots\varphi_{i_1+t-1}^{(in)}\s(i,\g,x_n+T)\s(n,-\frac{1}{\g},\frac{1}{\g}(x_i-T)),
\eee
where $\varphi_j^{(in)}$ are elementary automorphisms, according to (\ref{af4}). To (\ref{af18}) corresponds the sequence of elementary transformations  
\bes
(x_1,\ldots,x_{i-1},x_i,x_{i+1},\ldots,x_{n-1},x_n)\mapsto(x_1,\ldots,x_{i-1},x_n,x_{i+1},\ldots,x_{n-1},x_i)
\rightarrow \\
\te_1'\rightarrow\te_2'\rightarrow 
\ldots\rightarrow\te_{i_1+t-1}'
=\psi_1 \rightarrow \psi_2\rightarrow \te, 
\ees
where $\te_i'$ is obtained from $\te_i$ only by the permutation of the coordinates with numbers $i$ and $n$, 
and the transformation  
\bes
(x_1,\ldots,x_{i-1},x_i,x_{i+1},\ldots,x_{n-1},x_n)\mapsto(x_1,\ldots,x_{i-1},x_n,x_{i+1},\ldots,x_{n-1},x_i)
\ees
 is a composition of three elementary 
linear transformations. 

If in (\ref{af9}) we replace $\te$ by (\ref{af18}), 
then the exponent of (\ref{af10}) remains the same. But instead of  
$\phi$ we have $\psi_2$, which reduces the element $f_n$ of $\te$, and Lemma \ref{al6} gives the triviality of (\ref{af9}). $\Box$

\begin{lm}\label{al13}
If the elements $\overline{f_1},\ldots,\overline{f_{i-1}},\overline{f_{i+1}},\ldots,\overline{f_n}$ are free, then (\ref{af9}) is trivial.
\end{lm}
\Proof By Lemma \ref{al3} and (\ref{af16}), we have 
\bes
\overline{f_i} = \overline{C} = c(\overline{f_1},\ldots,\overline{f_{i-1}},\overline{f_{i+1}},\ldots,\overline{f_n}). 
\ees
By Lemma \ref{al9} we can assume that 
 $\overline{f_i}$ depends on $\overline{f_n}$. Consequently, 
$\deg f_n \leq \deg f_i$. 

Note that if the elements $\overline{f_1},\ldots,\overline{f_{i-1}},\overline{f_i},\overline{f_{i+1}},\ldots,\overline{f_{n-1}}$
 are not free, then it follows that  the elements $\overline{f_1},\ldots,\overline{f_{i-1}},\overline{f_{i+1}},\ldots,\overline{f_n}$
 are also not free, which contradicts the condition of the lemma.  
Consequently, the elements $\overline{f_1},\ldots,\overline{f_{i-1}},\overline{f_i},\overline{f_{i+1}},\ldots,\overline{f_{n-1}}$ are free. 
Then $\overline{B}\in \langle \overline{f_1},\ldots,\overline{f_{i-1}},\overline{f_i},\overline{f_{i+1}},\ldots,\overline{f_{n-1}} \rangle$. 
By Lemma \ref{al11} we can assume that $B$ contains $f_i$. 
Then $\deg\,f_i\leq \deg\,B\leq \deg\,f_n $, i.e.  $\deg\,f_i=\deg\,f_n$. 
 Hence  
\bes
\overline{C}=\overline{f_i}=\g \overline{f_n}+T(\overline{f_1},\ldots,\overline{f_{i-1}},\overline{f_{i+1}},\ldots,\overline{f_{n-1}}).
\ees
Lemma \ref{al12} gives the triviality of (\ref{af9}). $\Box$

\begin{lm}\label{al14}
Assume that there exists $r$ such that $r\neq i$, $1\leq i\leq n-1$, and 
\bes
\overline{f_r} \in <\overline{f_1},\ldots,\overline{f_{r-1}},\overline{f_{r+1}},\ldots,\overline{f_{i-1}},\overline{f_{i+1}},\ldots,\overline{f_n}>. 
\ees 
Then in (\ref{af11}) the automorphism $\phi$ can be replaced by an automorphism which reduces the element $f_r$ of $\te$. 
\end{lm}
\Proof 
Assume that 
\bes
\overline{f_r} = T(\overline{f_1},\ldots,\overline{f_{r-1}},\overline{f_{r+1}},\ldots,\overline{f_{i-1}},\overline{f_{i+1}},\ldots,\overline{f_n}) 
\ees
where $T\in F_{\mathfrak M}<X\setminus \{x_r,x_i\}>$, and put 
\bes
 g_r=f_r-T(f_1,\ldots,f_{r-1},f_{r+1},\ldots,f_{i-1},f_{i+1},\ldots,f_n). 
\ees

  By (\ref{af3}) we have  
\bes
\s_1=\s(r,1,-T)\s(i,1,c_1)\s(r,1,T)
\ees
for some $c_1\in F_{\mathfrak M}<X\setminus \{x_i\}>$.
After such replacement, instead of $\phi\rightarrow\te$ we obtain   
\bes
\phi\rightarrow \psi_1\rightarrow \psi_2\rightarrow \te,  
 \ees
where 
\bes
\psi_1=(f_1,\ldots,f_{r-1},g_r,f_{r+1},\ldots,f_{i-1},g_i,f_{i+1},\ldots,f_n), \\
\psi_2=(f_1,\ldots,f_{r-1},g_r,f_{r+1},\ldots,f_{i-1},f_i,f_{i+1},\ldots,f_n). 
\ees

Since $\deg\,\psi_1,\deg\,\psi_2<d=\deg\,\te$, the new sequence (\ref{af10}) has the same exponent.  
Then instead of $\phi$ in (\ref{af11}) we have $\psi_2$, which reduces the element $f_r$ of $\te$. $\Box$

{\bf Completion of the proof of Theorem \ref{at1}.} By Lemma \ref{al13}, we can assume that the elements 
\bes
\overline{f_1},\ldots,\overline{f_{i-1}},\overline{f_{i+1}},\ldots,\overline{f_n}
\ees
are not free. Then, according to Lemma \ref{al4}, there exists $j\neq i$ such that 
\bes
\overline{f_j} \in <\overline{f_1},\ldots,\overline{f_{j-1}},\overline{f_{j+1}},\ldots,\overline{f_{i-1}},\overline{f_{i+1}},\ldots,\overline{f_n}>. 
\ees
By Lemma \ref{al10}, we can assume that $j\neq n$, i.e. $j\leq n-1$. According to Lemma \ref{al8}, 
we can also assume that $\overline{f_j}$ depends on $\overline{f_n}$. Consequently, $\deg\,f_j\geq \deg\,f_n$. 
If $\deg\,f_j=\deg\,f_n$, then from this we can easily obtain that 
$\overline{f_n} \in <\overline{f_1},\ldots,\overline{f_{i-1}},\overline{f_{i+1}},\ldots,\overline{f_{n-1}}>$, 
and Lemma \ref{al10} gives the triviality of (\ref{af9}). Thus, it can be assumed that $\deg\,f_j> \deg\,f_n$. Moreover, 
by Lemma \ref{al14}, we may assume that $\phi$ reduces the element $f_j$ of $\te$. Interchanging $f_i$ and $f_j$, 
from now we can assume without of generality that $\deg\,f_i>\deg\,f_n$ and 
\bee\label{af19}
\overline{f_i} \in <\overline{f_1},\ldots,\overline{f_{i-1}},\overline{f_{i+1}},\ldots,\overline{f_n}>. 
\eee
Suppose that the elements $\overline{f_1},\overline{f_2},\ldots,\overline{f_{n-1}}$ are free. By Lemma \ref{al3}, 
$\overline{B}\in <\overline{f_1},\overline{f_2},\ldots,\overline{f_{n-1}}>$.  If $B$ depends on $f_i$, then $\deg\,B \geq \deg\,f_i > \deg\,f_n$, which contradicts 
(\ref{af15}). If $B$ does not depend on $f_i$, then Lemma \ref{al11} gives the triviality of (\ref{af9}). 

If the elements $\overline{f_1},\overline{f_2},\ldots,\overline{f_{n-1}}$ are not free, then there exists $r$ such that 
\bee\label{af20}
\overline{f_r} \in <\overline{f_1},\ldots,\overline{f_{r-1}},\overline{f_{r+1}},\ldots,\overline{f_{n-1}}>. 
\eee
By Lemma \ref{al9}, we can take $r\neq i$. If $\overline{f_r}$ does not depend on $\overline{f_i}$, 
then Lemma \ref{al8} gives the triviality of (\ref{af9}). Assume that $\overline{f_r}$  depends on $\overline{f_i}$. Then, $\deg\,f_r\geq \deg\,f_i$.  
If $\deg\,f_r= \deg\,f_i$, then from this we can obtain that  $\overline{f_i} \in <\overline{f_1},\ldots,\overline{f_{i-1}},\overline{f_{i+1}},\ldots,\overline{f_{n-1}}>$, 
and Lemma \ref{al9} gives the triviality of (\ref{af9}). So, we can assume that $\deg\,f_r > \deg\,f_i$. 
Then (\ref{af19}) gives that 
\bes
\overline{f_i} \in <\overline{f_1},\ldots,\overline{f_{r-1}},\overline{f_{r+1}},\ldots,\overline{f_{i-1}},\overline{f_{i+1}},\ldots,\overline{f_n}>.  
\ees
Consequently, 
\bes
\overline{f_r} \in <\overline{f_1},\ldots,\overline{f_{r-1}},\overline{f_{r+1}},\ldots,\overline{f_{i-1}},\overline{f_{i+1}},\ldots,\overline{f_n}>.  
\ees
By Lemma \ref{al14}, we can assume that $\phi$ reduces the element $f_r$ of $\te$. Then, (\ref{af20}) and Lemma \ref{al9} gives the triviality of (\ref{af9}). 

This completes the proof of Theorem \ref{at1}. $\Box$

\begin{center}
{\bf\large Acknowledgments} 
\end{center}

\hspace*{\parindent}

I am  grateful to Max-Planck Institute f\"ur Mathematik for hospitality and exellent working conditions. 
I am also grateful to V.\,Shpilrain, M.\,Zaidenberg for very helpful discussions and comments.

\bigskip

\hspace*{\parindent}

\end{document}